\documentclass[11pt,a4]{article}

\usepackage{graphicx}

\def\be{\begin{equation}}
\def\ee{\end{equation}}
\newcommand{\kk}[2]{\frac{#1}{#2}}
\newcommand{\ff}[1]{{\bf  #1}}
\def\a{\alpha}
\def\b{\beta}
\def\x{\ff{x}}
\def\stepa{\indent \quad \quad}
\def\stepb{\indent \qquad \qquad}
\newcommand{\vcode}[4]{\begin{figure}
\begin{center}
\begin{minipage}[c]{#1\textwidth}
{{\small #2 \hrule   %
{\it #3} \\[-5pt] \hrule }}
\end{minipage}
\caption{#4}
\end{center}   \end{figure}} 

\begin{document}

\title{A New Metaheuristic Bat-Inspired Algorithm}

\author{Xin-She Yang \\
Department of Engineering, University of Cambridge, \\
Trumpington Street, Cambridge CB2 1PZ, UK \\
}

\date{}

\maketitle

\abstract{
Metaheuristic algorithms such as particle swarm optimization, firefly algorithm
and harmony search are now becoming powerful methods for solving many tough optimization problems.
In this paper, we propose a new metaheuristic method, the Bat Algorithm,
based on the echolocation behaviour of bats. We also intend to
combine the advantages of existing algorithms into the new bat algorithm.
After a detailed formulation and explanation of its implementation,
we will then compare the proposed algorithm with other existing algorithms,
including genetic algorithms and particle swarm optimization.
Simulations show that the proposed algorithm seems
much superior to other algorithms, and further studies are also discussed. \\
}

\noindent {\bf Citation detail:}

X.-S. Yang,  A New Metaheuristic Bat-Inspired Algorithm, in: {\it Nature Inspired Cooperative Strategies for Optimization (NISCO 2010)} (Eds. J. R. Gonzalez et al.), Studies in Computational Intelligence,
Springer Berlin, {\bf 284}, Springer, 65-74 (2010).

\section{Introduction}

Metaheuristic algorithms such as particle swarm optimization and simulated annealing
are now becoming powerful methods for solving many tough optimization problems [3-7,11].
The vast majority of heuristic and metaheuristic algorithms have been derived from the
behaviour of biological systems and/or physical systems in nature. For example, particle swarm
optimization was developed based on the swarm behaviour of birds and fish \cite{Kennedy,Kennedy2},
while simulated annealing was based on the annealing process of metals \cite{Kirk}.

New algorithms are also emerging recently, including harmony search
and the firefly algorithm. The former was inspired by the improvising process of composing
a piece of music \cite{Geem}, while the latter was formulated based on the flashing
behaviour of fireflies \cite{Yang}. Each of these algorithms has certain advantages
and disadvantages. For example, simulating annealing can almost guarantee to find the
optimal solution if the cooling process is slow enough and the simulation is running long
enough; however, the fine adjustment in parameters does affect the convergence rate of
the optimization process. A natural question is whether it is possible to combine major advantages
of these algorithms and try to develop a potentially better algorithm? This paper is such
an attempt to address this issue.

In this paper, we intend to propose a new metaheuristic method, namely, the Bat Algorithm (BA),
based on the echolocation behaviour of bats. The capability of echolocation of
microbats is fascinating as these bats can find their prey and discriminate different
types of insects even in complete darkness. We will first formulate the bat algorithm
by idealizing the echolocation behaviour of bats. We then describe how it works
and make comparison with other existing algorithms. Finally, we will discuss some
implications for further studies.

\section{Echolocation of bats}

\subsection{Behaviour of microbats}

Bats are fascinating animals. They are the only mammals with wings and
they also have advanced capability of echolocation.
It is estimated that there are about 996
different species which account for up to 20\% of all mammal
species \cite{Alt,Colin}. Their size ranges from the tiny bumblebee bat
(of about 1.5 to 2g) to the giant bats with wingspan of about 2 m and weight
up to about 1 kg. Microbats typically have forearm length of
about 2.2 to 11cm. Most bats uses echolocation to a certain degree;
among all the species, microbats are a famous example as microbats use
echolocation extensively while megabats do not \cite{Bat,Bat2}.

Most microbats are insectivores.
Microbats use a type of sonar, called, echolocation, to detect prey, avoid
obstacles, and locate their roosting crevices in the dark. These bats emit
a very loud sound pulse and listen for the echo that bounces back from
the surrounding objects. Their pulses vary in properties and can be
correlated with their hunting strategies, depending on the species.
Most bats use short, frequency-modulated signals to sweep through about an octave,
while others more often use constant-frequency signals for echolocation. Their signal
bandwidth varies depends on the species, and often increased by using more
harmonics.

\subsection{Acoustics of Echolocation}

Though each pulse only lasts a few thousandths of
a second (up to about 8 to 10 ms), however, it has a constant frequency which
is usually in the region of 25kHz to 150 kHz. The typical range of frequencies
for most bat species are in the region between 25kHz and 100kHz,
though some species can emit higher frequencies up to 150 kHz.
Each ultrasonic burst may last typically 5 to 20 ms, and microbats emit
about 10 to 20 such sound bursts every second. When hunting for prey, the rate
of pulse emission can be sped up to
about 200 pulses per  second when they fly near their prey.
Such short sound bursts imply the fantastic ability of the signal
processing power of bats. In fact, studies shows the integration
time of the bat ear is typically about 300 to 400 $\mu$s.

As the speed of sound in air is typically $v=340$ m/s, the wavelength
$\lambda$ of the ultrasonic sound bursts with a constant frequency $f$ is given by
\be \lambda=\frac{v}{f}, \ee
which is in the range of 2mm to 14mm for the typical frequency range
from 25kHz to 150 kHz. Such wavelengths are in the
same order of their prey sizes.

Amazingly, the emitted pulse could be
as loud as 110 dB, and, fortunately, they are in the ultrasonic
region.  The loudness also varies from the loudest when searching
for prey and to a quieter base when homing towards the prey.
The travelling range of such short
pulses are typically a few metres, depending on the actual frequencies \cite{Bat}.
Microbats can manage to avoid obstacles as small as thin human hairs.

Studies show that microbats use the time delay from the emission
and detection of the echo, the time difference between their two
ears, and the loudness variations of the echoes to build up three
dimensional scenario of the surrounding. They can detect the distance
and orientation of the target, the type of prey, and even the moving
speed of the prey such as small insects. Indeed, studies suggested
that bats seem to be able to discriminate targets by the variations of the
Doppler effect induced by the wing-flutter rates of the target
insects \cite{Alt}.

Obviously, some bats have good eyesight, and most bats also have
very sensitive smell sense. In reality, they will use all the senses
as a combination to maximize the efficient detection of prey and
smooth navigation. However, here we are only interested in the echolocation
and the associated behaviour.

Such echolocation behaviour of microbats
can be formulated in such a way that it can be associated with
the objective function to be optimized, and this make it possible
to formulate new optimization algorithms. In the rest of this paper,
we will first outline the basic formulation of the
Bat Algorithm (BA) and then discuss
the implementation  and comparison in detail.

\section{Bat Algorithm}

If we idealize some of the echolocation characteristics of microbats,
we can develop various bat-inspired algorithms or bat algorithms.
For simplicity, we now use the following approximate or idealized
rules:

\begin{itemize}

\item[1.] All bats use echolocation to sense distance, and
they also `know' the difference between food/prey and background barriers
in some magical way;

\item[2.] Bats fly randomly with velocity $\ff{v}_i$
at position $\x_i$ with a fixed frequency $f_{\min}$, varying wavelength $\lambda$
and loudness $A_0$ to search for prey. They can automatically adjust
the wavelength (or frequency) of their
emitted pulses and adjust the rate of pulse emission $r \in [0,1]$,
depending on the proximity of their target;

\item[3.] Although the loudness can vary in many ways, we assume that
the loudness varies from a large (positive) $A_0$ to a minimum
constant value $A_{\min}$.

\end{itemize}

Another obvious simplification is that no ray tracing is used in
estimating the time delay and three dimensional topography. Though this
might be a good feature for the application in computational geometry,
however, we will not use this as it is more computationally extensive in
multidimensional cases.

In addition to these simplified assumptions, we also use the following
approximations, for simplicity. In general the frequency $f$  in
a range $[f_{\min}, f_{\max}]$ corresponds to a range
of wavelengths $[\lambda_{\min}, \lambda_{\max}]$. For example a frequency
range of [$20$kHz, $500$kHz] corresponds to a range of wavelengths
from $0.7$mm to $17$mm.

For a given problem, we can also  use any wavelength  for the ease of implementation.
In the actual implementation, we can
adjust the range by adjusting the wavelengths (or frequencies), and
the detectable range (or the largest wavelength) should be chosen such
that it is comparable to the size of the domain of interest, and then toning
down to smaller ranges.
Furthermore, we do not necessarily have to use the wavelengths themselves,
instead, we can also vary the frequency while fixing the wavelength $\lambda$.
This is because $\lambda$ and $f$ are related due to the fact $\lambda f$ is constant.
We will use this later approach in our implementation.

For simplicity, we can assume
$f \in [0, f_{\max}]$. We know that higher frequencies
have short wavelengths and travel a shorter distance. For bats, the
typical ranges are a few metres.
The rate of pulse can simply be in the range of $[0,1]$ where $0$ means
no pulses at all, and $1$ means the maximum rate of pulse emission.

Based on these approximations and idealization, the basic steps of the Bat Algorithm (BA)
can be summarized as the pseudo code shown in Fig. \ref{fa-fig-100}.

\vcode{0.9}{{\sf Bat Algorithm}} {
\stepa Objective function $f(\x), \quad \x=(x_1, ..., x_d)^T$ \\
\stepa Initialize  the bat population $\x_i \; (i=1,2,...,n)$ and  $\ff{v}_i$\\
\stepa Define pulse frequency $f_i$ at $\x_i$  \\
\stepa Initialize pulse rates $r_i$ and the loudness $A_i$ \\
\stepa {\bf while} ($t<$Max number of iterations) \\
\stepa Generate new solutions by adjusting frequency,  \\
\stepa and updating velocities and locations/solutions
[equations (\ref{f-equ-150}) to (\ref{f-equ-250})]  \\
\stepb {\bf if } (rand $>r_i$)    \\
\stepb Select a solution among the best solutions \\
\stepb Generate a local solution around the selected best solution \\
\stepb {\bf end if} \\
\stepb Generate a new solution by flying randomly \\
\stepb {\bf if  } (rand $<A_i \;\; \&  \;\; f(\x_i)<f(\x_*)$) \\
\stepb Accept the new solutions \\
\stepb Increase $r_i$ and reduce $A_i$ \\
\stepb {\bf end if} \\
\stepa Rank the bats and find the current best $\x_*$ \\
\stepa {\bf end while} \\
\stepa Postprocess results and visualization }
{Pseudo code of the bat algorithm (BA). \label{fa-fig-100} }

\subsection{Movement of Virtual Bats}

In simulations, we use virtual bats naturally.
We have to define the rules how their
positions $\x_i$ and velocities $\ff{v}_i$ in a $d$-dimensional search space
are updated.
The new solutions $\x_i^{t}$ and velocities $\ff{v}_i^{t}$ at time step
$t$ are given by
\be f_i =f_{\min} + (f_{\max}-f_{\min}) \beta, \label{f-equ-150} \ee
\be \ff{v}_i^{t} = \ff{v}_i^{t-1} +  (\x_i^t - \x_*) f_i , \ee
\be \x_i^{t}=\x_i^{t-1} + \ff{v}_i^t,  \label{f-equ-250} \ee
where $\beta \in [0,1]$ is a random vector drawn from a uniform distribution.
Here $\x_*$ is the current global best
location (solution) which is located after comparing all
the solutions among all the $n$ bats. As the product $\lambda_i f_i$ is the
velocity increment, we can use either $f_i$ (or $\lambda_i$ ) to adjust the velocity change
while fixing the other factor $\lambda_i$ (or $f_i$), depending on the
type of the problem of interest.
In our implementation, we will use $f_{\min}=0$ and $f_{\max}=100$, depending the
domain size of the problem of interest. Initially, each bat is randomly assigned
a frequency which is drawn uniformly from $[f_{\min}, f_{\max}]$.

For the local search part, once a solution is selected among the current best solutions,
a new solution for each bat is generated locally using random walk
\be \x_{\rm new}=\x_{\rm old} + \epsilon  A^{t}, \ee
where $\epsilon \in [-1,1]$ is a random number, while $A^{t}=<\!\!A_i^{t}\!>$
is the average loudness of all the bats at this time step.

The update of the velocities
and positions of bats have some similarity to the procedure in  the standard
particle swarm optimization \cite{Kennedy} as $f_i$ essentially
controls the pace and range of the movement of the swarming particles.
To a degree, BA can be considered
as a balanced combination of the standard particle swarm optimization and
the intensive local search controlled by the loudness and pulse rate.

\subsection{Loudness and Pulse Emission}

Furthermore, the loudness $A_i$ and the rate $r_i$ of pulse emission have to be
updated accordingly as the iterations proceed. As the loudness usually
decreases once a bat has found its prey, while the rate of pulse emission
increases, the loudness can be chosen as any value of convenience.
For example, we can use
$A_0=100$ and $A_{\min}=1$. For simplicity, we can also use
$A_0=1$ and $A_{\min}=0$, assuming $A_{\min}=0$ means that a bat
has just found the prey and temporarily stop emitting any sound.
Now we have
\be A_i^{t+1}=\alpha A_{i}^{t}, \;\;\;\;\; r_i^{t+1}= r_i^0 [1-\exp(-\gamma t)],
\label{rate-equ-50} \ee
where $\alpha$ and $\gamma$ are constants. In fact, $\alpha$ is similar
to the cooling factor of a cooling schedule in the simulated annealing \cite{Kirk}.
For any $0<\alpha<1$ and $\gamma>0$, we have
\be A_i^t \rightarrow 0, \;\;\; r_i^t \rightarrow r_i^0, \;\;\textrm{as} \;\;
t \rightarrow \infty. \ee In the simplicity case, we can use $\a=\gamma$,
and we have used $\a=\gamma=0.9$ in our simulations. The choice of parameters
requires some experimenting.
Initially, each bat should have different values of loudness and pulse emission rate,
and this can be achieved by randomization. For example, the initial loudness $A_i^{0}$
can typically be $[1, 2]$, while the initial emission rate $r_i^{0}$ can be around
zero, or any value $r_i^0 \in [0,1]$ if using (\ref{rate-equ-50}).
Their loudness and emission rates
will be updated only if the new solutions are improved, which means that these bats
are moving towards the optimal solution.

\section{Validation and Comparison}

From the pseudo code, it is relatively straightforward to implement
the Bat Algorithm in any programming language. For the ease of visualization,
we have implemented it using Matlab for various test functions.

\subsection{Benchmark Functions}

There are many standard test functions for validating new algorithms.
In the current benchmark validation, we have chosen the well-known Rosenbrock's function
\be f(\x)=\sum_{i=1}^{d-1} (1-x_i^2)^2+100 (x_{i+1}-x_i^2)^2, \;\; -2.048 \le x_i \le 2.048, \ee
and the eggcrate function
\be g(x,y)=x^2+y^2+25 (\sin^2 x + \sin^2 y), \;\;\; (x,y) \in [-2 \pi, 2 \pi] \times
[-2 \pi, 2 \pi]. \ee
We know that $f(\x)$ has a global minimum $f_{\min}=0$ at $(1,1)$ in 2D, while
$g(x,y)$ has a global minimum $g_{\min}=0$ at $(0,0)$.
De Jong's standard sphere function
\be h(\x)=\sum_{i=1}^d x_i^2, \;\;\; -10 \le x_i \le 10, \ee
has also been used. Its minimum is $h_{\min}=0$ at $(0,0,...,0)$ for any $d \ge 3$.

\begin{figure}
 \centerline{\includegraphics[height=3.5in,width=5in]{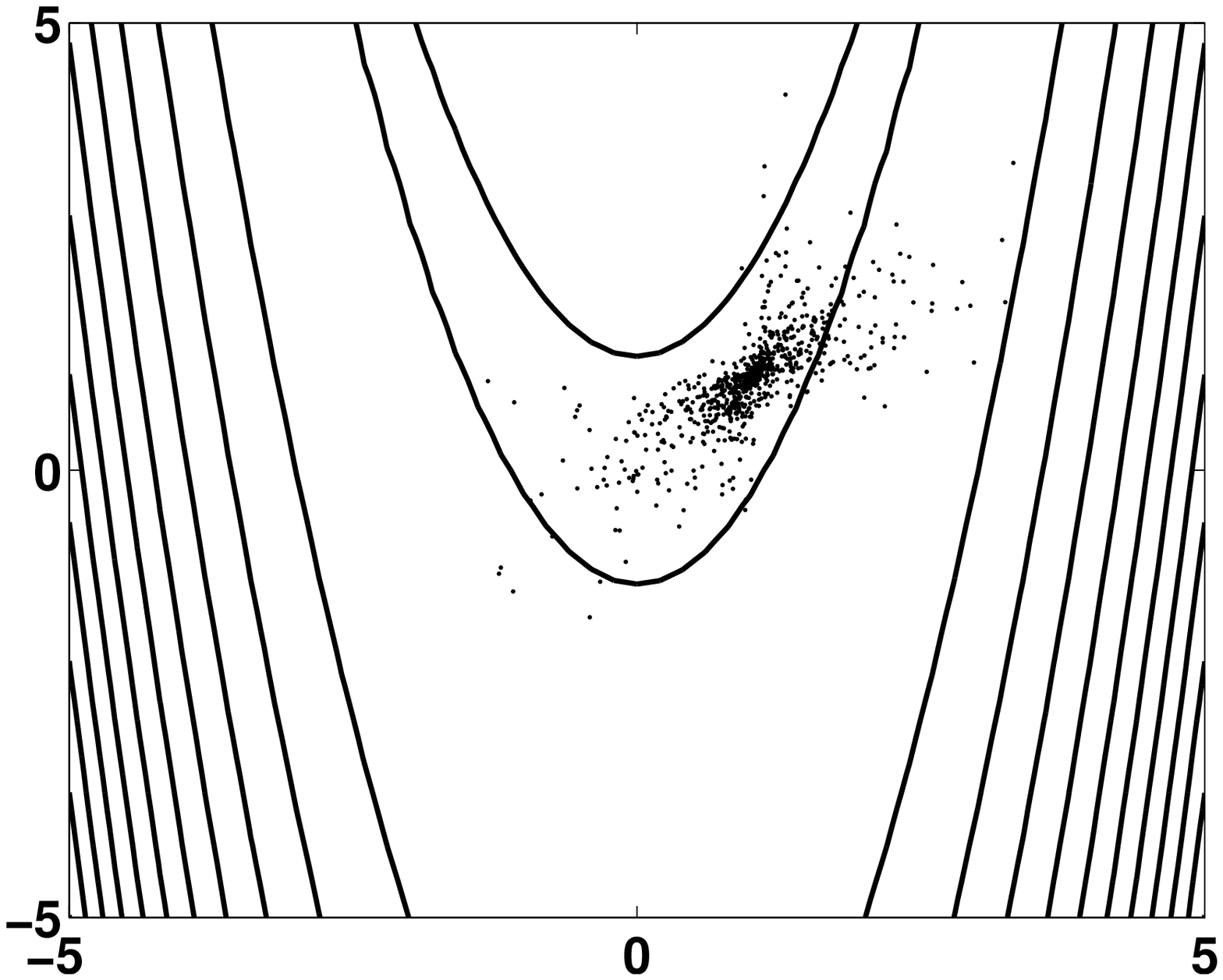}}
\vspace{-5pt}
\caption{The paths of 25 virtual bats during 20 consecutive iterations.
They converge into $(1,1)$.  \label{ba-fig-1} }
\end{figure}

In addition, we have also used other standard test functions for numerical global optimization \cite{Liang}
such as Ackley's function
\be s(\x)=20+e -20 \exp \Big[-0.2 \sqrt{\frac{1}{d} \sum_{i=1}^d x_i^2} \Big]
-\exp[\frac{1}{d} \sum_{i=1}^d \cos (2 \pi x_i)], \ee
where $-30 \le x_i \le 30$.  It has the global minimum $s_{\min}=0$ at $(0,0,...,0)$.

Michalewicz's test function
\be f(\x)=-\sum_{i=1}^d \sin (x_i) \Big[ \sin (\kk{i x_i^2}{\pi}) \Big]^{2m}, \;\; (m=10), \ee
has $d!$ local optima in the the domain $0 \le x_i \le \pi$ where $i=1,2,...,d$.
The global minimum is $f_* \approx -1.801$ for $d=2$, while $f_* \approx -4.6877$ for $d=5$.

\begin{figure}
 \centerline{\includegraphics[height=2.5in,width=3in]{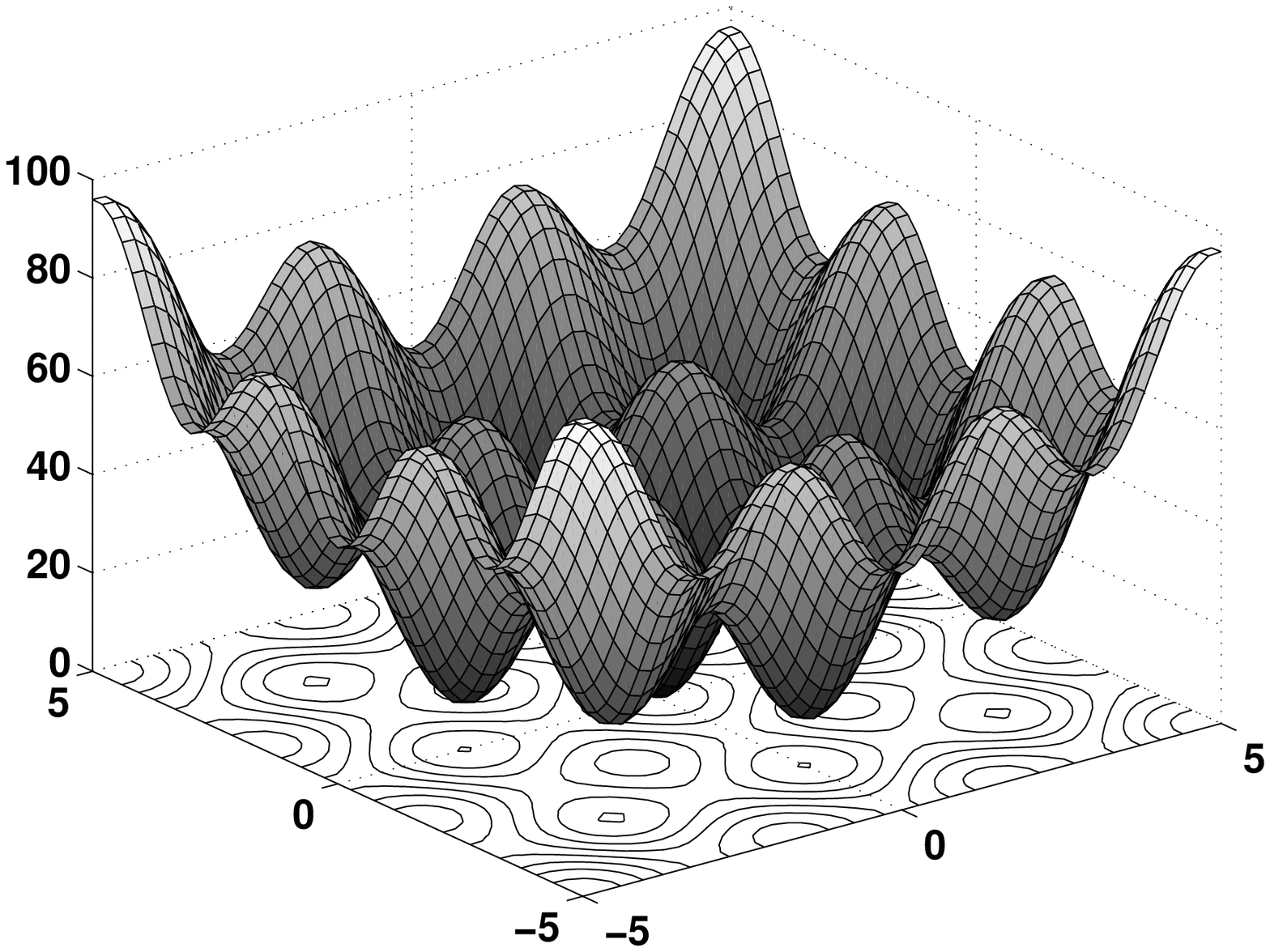}
 \includegraphics[height=2.5in,width=3in]{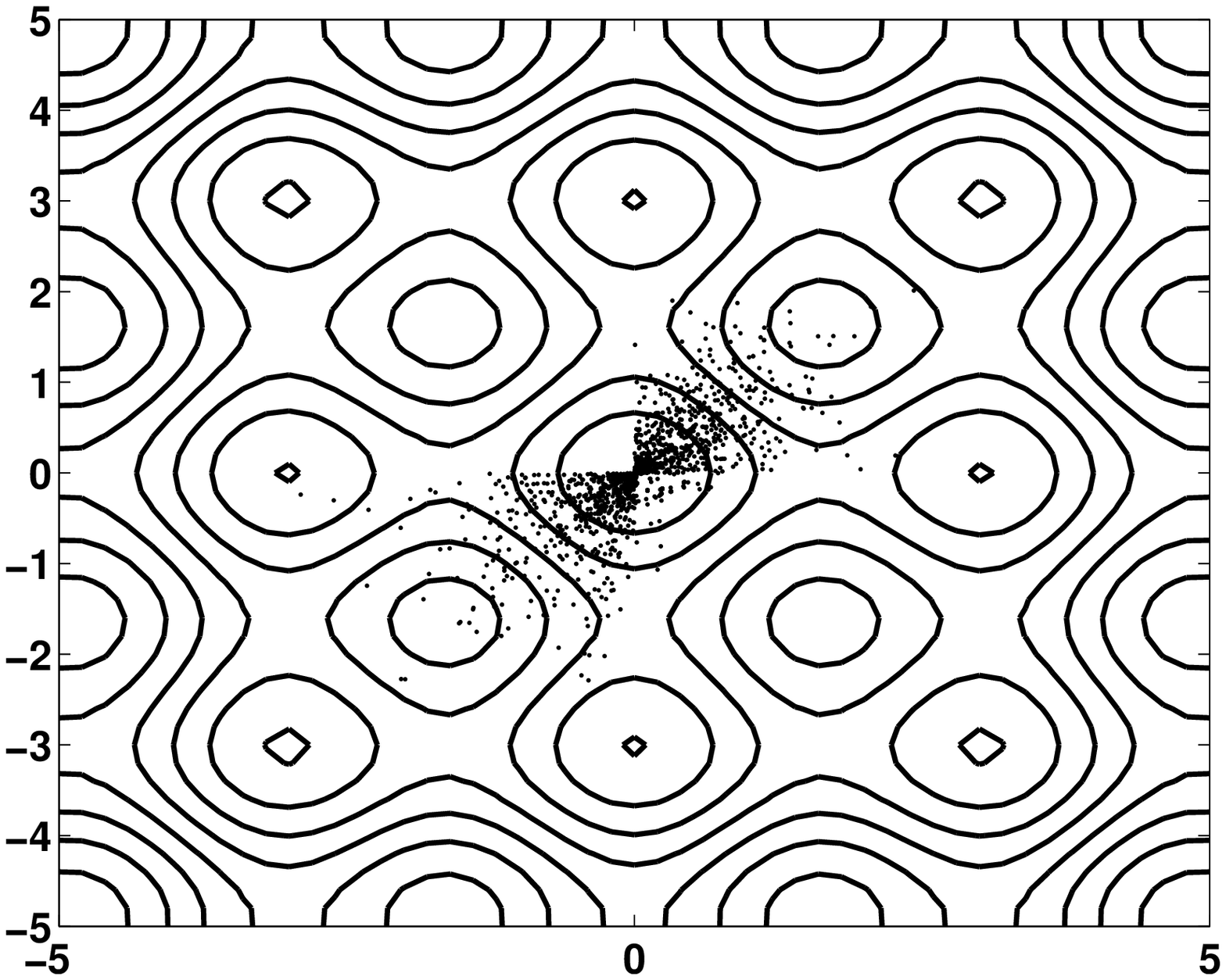}}
\vspace{-5pt}
\caption{The eggcrate function (left) and
the locations of 40 bats in the last ten iterations (right).  \label{ba-fig-2} }
\end{figure}

In our implementation, we use $n=25$ to $50$ virtual bats, and $\alpha=0.9$. For Rosenbrock's 2-D banana function, the paths
of 25 virtual bats during the consecutive 20 time steps are shown in Fig. \ref{ba-fig-1}
where we can see that the bats converge at the global optimum $(1,1)$.
For the multimodal eggcrate function, a snapshot of the last 10 iterations is shown
in Fig. \ref{ba-fig-2}. Again, all bats move towards the global best $(0,0)$.

\subsection{Comparison with Other Algorithms}

In order to compare the performance of the new algorithm, we have tested it against
other heuristic algorithms, including genetic algorithms (GA) \cite{Gold,Holland,Mitchell},
and particle swarm optimization (PSO) \cite{Kennedy,Kennedy2}. There are many variants
of PSO, and some variants such as the mean PSO could perform better than
the standard PSO \cite{Deep}; however, the standard PSO is by far the most popularly used.
Therefore, we will also use the standard PSO in our comparison.

There are many ways to carry out the comparison of algorithm performance, and two obvious
approaches are: to compare the numbers of function evaluations for a given tolerance or accuracy, or
to compare their accuracies for a fixed number of function evaluations. Here we will use
the first approach. In our simulations, we use a fixed tolerance $\epsilon \le 10^{-5}$, and
we run each algorithm for $100$ times so that we can do meaningful statistical analysis.

For genetic algorithms, we have used the standard version with no elitism
with the mutation probability of $p_m=0.05$ and crossover probability of $0.95$.
For particle swarm optimization, we have also used the standard version with
learning parameters $\a=\b=2$ and the inertia function $I=1$ \cite{Kennedy,Kennedy2}.
The simulations have been carried out using Matlab on a standard 3GHz desktop computer.
Each run with about 10,000 function evaluations typically takes less
than 5 seconds. Furthermore, we have tried to use different population sizes
from $n=10$ to $250$, and we found that for most problems, $n=15$ to $50$ is
sufficient. Therefore, we use a fixed population $n=40$ for all simulations.
Table 1 shows the number of function evaluations in the form
of mean $\pm$ the standard deviation (success rate of the algorithm in finding the
global optima).

\begin{table}[ht]

\caption{Comparison of BA with GA, and PSO.}
\centering

\begin{tabular}{ccccc}
\hline \hline
Functions/Algorithms & GA & PSO & BA \\
\hline

 Multiple peaks & $52124 \pm 3277 (98\%)$  & $3719 \pm  205 (97\%)$ & $1152 \pm 245 (100 \%)$ \\

 Michalewicz's ($d\!\!=\!\!16$)  & $89325 \pm 7914 (95 \%)$  & $6922 \pm 537 (98\%)$  & $4752 \pm 753 (100\%)$ \\

 Rosenbrock's ($d\!\!=\!\!16$) & $55723 \pm 8901 (90\%)$ & $32756 \pm 5325 (98\%)$ & $7923 \pm 3293 (100\%) $ \\

 De Jong's ($d\!\!=\!\!256$) & $25412 \pm 1237 (100\%)$ & $17040 \pm 1123 (100\%)$ & $5273 \pm 490 (100\%)$\\
 Schwefel's ($d\!\!=\!\!128$) & $227329 \pm 7572 (95\%)$ & $14522 \pm 1275 (97\%)$ & $8929 \pm 729 (99\%)$ \\

 Ackley's ($d\!\!=\!\!128$) & $32720 \pm 3327 (90\%)$ & $23407 \pm 4325 (92\%)$ & $6933 \pm 2317 (100\%)$ \\

 Rastrigin's & $110523 \pm 5199 (77 \%)$ & $79491 \pm 3715 (90\%)$ & $12573 \pm 3372 (100\%)$ \\

 Easom's & $19239 \pm 3307 (92\%)$ & $17273 \pm 2929 (90\%)$ & $7532 \pm 1702 (99\%)$ \\

 Griewangk's & $70925 \pm 7652 (90\%)$ & $55970 \pm 4223 (92\%)$ & $9792 \pm 4732 (100\%)$ \\

 Shubert's (18 minima) & $54077 \pm 4997 (89\%)$ & $23992 \pm 3755 (92\%)$ & $11925 \pm 4049 (100\%)$ \\

\hline
\end{tabular}
\end{table}

From Table 1, we can see that PSO performs much better than genetic algorithms,
while the Bat Algorithm is much superior to other
algorithms in terms of accuracy and efficiency. This is no surprising as the aim of developing the new
algorithm was to try to use the advantages of existing algorithms and other interesting feature
inspired by the fantastic behaviour of echolocation of microbats.

If we replace the variations of the frequency $f_i$ by a
random parameter and setting $A_i=0$ and $r_i=1$,
the bat algorithm essentially becomes the standard Particle Swarm Optimization (PSO).
Similarly, if we do not use the velocities, but we use fixed loudness and rate:
$A_i$ and $r_i$. For example, $A_i=r_i=0.7$, this algorithm is virtually
reduced to a simple Harmony Search (HS) as the frequency/wavelength change is
essentially the pitch adjustment, while the rate of pulse emission is similar to the harmonic
acceptance rate (here with a twist) in the harmony search algorithm \cite{Geem,Yang2}.
The current studies implies that the proposed new algorithm is potentially more powerful
and thus should be investigated further in many applications of engineering and industrial
optimization problems.

\section{Discussions}

In this paper, we have successfully formulated a new Bat Algorithm for
continuous constrained optimization problems.
From the formulation of the Bat Algorithm and its implementation and comparison,
we can see that it is a very promising algorithm. It is potentially
more powerful than particle swarm optimization and genetic algorithms
as well as Harmony Search. The primary reason is that BA uses
a good combination of major advantages of these algorithms
in some way. Moreover, PSO and harmony search are the special cases
of the Bat Algorithm under appropriate simplifications.

In addition, the fine adjustment of the parameters $\alpha$ and $\gamma$
can affect the convergence rate of the bat algorithm. In fact, parameter
$\alpha$ acts in a similar role as the cooling schedule in the simulated
annealing. Though the implementation is more complicated than many other
metaheuristic algorithms; however, it does show that it utilizes a balanced
combination of the advantages of existing successful algorithms with
innovative feature based on the echolocation behaviour of bats.
New solutions are generated by adjusting frequencies, loudness and
pulse emission rates, while the proposed solution is accepted or not
depends on the quality of the solutions controlled or characterized
by loudness and pulse rate which are in turn related to the closeness
or the fitness of the locations/solution to the global optimal solution.

The exciting results suggest that more studies will
highly be needed to carry out the sensitivity analysis,
to analyze the rate of algorithm convergence,
and to improve the convergence rate even further.
More extensive comparison studies with a more wide range of existing
algorithms using much tough test functions in
higher dimensions will pose more challenges to the algorithms, and
thus such comparisons will potentially reveal the virtues and weakness
of all the algorithms of interest.

An interesting extension will be to use different schemes
of wavelength or frequency variations
instead of the current linear implementation. In addition, the rates of pulse emission
and loudness can also be varied in a more sophisticated manner.
Another extension for discrete problems is to use the time delay
between pulse emission and the echo bounced back. For example,
in the travelling salesman problem, the distance between two adjacent
nodes/cities can easily be coded as time delay.
As microbats use time difference between their two ears to obtain three-dimensional
information, they can identify the type of prey and the velocity of a flying insect.
Therefore, a further natural extension to the current bat algorithm
would be to use the directional echolocation and Doppler effect,
which may lead to even more interesting variants and new algorithms.

\end{document}